\documentclass{IEEEtran}
\usepackage{graphicx}
\usepackage{amsmath}
\usepackage{amsfonts}
\usepackage{subfigure}
\usepackage{stmaryrd}
\usepackage{multirow, multicol, booktabs}
\usepackage{hyperref}
\usepackage{siunitx}
\DeclareMathOperator{\Div}{div}

\begin{document}
\title{Optimal Control of Hybrid Systems via Measure Relaxations}
\author{Etienne Buehrle$^{1}$, {\"O}mer~{\c{S}}ahin~Ta{\c{s}}$^{2}$, and Christoph Stiller$^{1}$
\thanks{$^{1}$Karlsruhe Institute of Technology, 76131 Karlsruhe, Germany
        {\tt\small \{etienne.buehrle, stiller\}@kit.edu}}
\thanks{$^{2}$FZI Research Center Information Technology, 76131 Karlsruhe, Germany
        {\tt\small tas@fzi.de}}}

\maketitle

\begin{abstract}
We propose an approach to trajectory optimization for piecewise polynomial systems based on the recently proposed graphs of convex sets framework. We instantiate the framework with a convex relaxation of optimal control based on occupation measures, resulting in a convex optimization problem resembling the discrete shortest-paths linear program that can be solved efficiently to global optimality. While this approach inherits the limitations of semidefinite programming, scalability to large numbers of discrete modes improves compared to the NP-hard mixed-integer formulation. We use this to plan trajectories under temporal logic specifications, comparing the computed cost lower bound to a nonconvex optimization approach with fixed mode sequence. In our numerical experiments, we find that this bound is typically in the vicinity of the nonconvex solution, while the runtime speedup is significant compared to the often intractable mixed-integer formulation. Our implementation is available at \url{https://github.com/ebuehrle/hpoc}.
\end{abstract}


\section{Introduction}
Systems integrating discrete reasoning and continuous control appear in a variety of contexts ranging from biological to robotic systems \cite{bemporad1999control, garrett2021integrated, tas2024decision, jamgochian2023shail}. A unified description of such mixed logical-dynamical systems is given by the class of piecewise polynomial systems \cite{bemporad2004efficient, sontag1981nonlinear}, which we consider in the following.

Multiple approaches to optimal control of piecewise-affine systems have been proposed. One line of work includes model predictive mixed-integer optimization \cite{bemporad1999control, bemporad2000explicit}, set-based methods \cite{althoff2021set, lygeros2004reachability}, and abstraction-based methods \cite{rungger2016scots, khaled2019pfaces, roy2011pessoa}. While the latter two approaches have proven challenging to scale to large state and action spaces owing to the curse of dimensionality, mixed-integer optimization is usually limited to systems with few discrete modes and short control horizons due to the combinatorial explosion of mode sequences.

Sampling-based approaches, although not subject to the curse of dimensionality a priori, are typically only probabilistically complete and may fail in large state-action spaces \cite{karaman2011sampling, lee2018monte}, or rely on non-convex optimization problems lacking convergence guarantees \cite{bertsekas2012dynamic}.

Prior work has considered convex relaxations of polynomial optimal control problems in discrete time \cite{han2019controller} and continuous time \cite{lasserre2008nonlinear, zhao2017control}. An extension to hybrid systems was proposed in \cite{zhao2019optimal}, however assuming identical input spaces in each mode.

The recently proposed graphs of convex sets \cite{marcucci2024graphs} have shown that a class of mixed-integer optimization problems on graphs can be reformulated as mixed-integer convex programs with tight convex relaxations.

We instantiate the graphs of convex sets framework with a convex relaxation of optimal control based on state-action occupation measures \cite{lasserre2008nonlinear}, yielding a compact convex formulation that is amenable to an efficient numerical solution.

\section{Background}
\label{sec:background}
We briefly recall the weak formulation of optimal control \cite{lasserre2008nonlinear} and the graphs of convex sets framework \cite{marcucci2024shortest}.

\subsection{Hybrid Optimal Control Problem}
\label{sec:ocp}
A piecewise polynomial system is modeled as a graph of modes $Q$ and transitions $\Sigma\subseteq Q \times Q$. Each mode $q\in Q$ is associated with a state space $X_q$, an input space $U_q$, system dynamics $f_q$, and a cost function $c_q$. This work considers the optimal control problem for piecewise polynomial systems, which seeks to determine cost-optimal state, input, and mode trajectories $x(t), u(t), q(t)$ and may be stated as
\begin{equation}\label{eq:pwa-ocp}\begin{aligned}
&\underset{x,u,q}{\text{minimize}} \quad &&\int_0^T c_q(x,u)dt \\
&\text{subject to} \quad &&\dot x = f_q(x,u) \\
&&&x \in X_q, u \in U_q
\end{aligned}\end{equation}
for a potentially variable control horizon $T$.

\subsection{Weak Optimal Control}
\label{sec:optimalcontrol}
For a single discrete mode, problem \eqref{eq:pwa-ocp} can be expressed as a minimization of a duality pairing $\langle c, \mu \rangle = \int_{X \times U} c\,d\mu$ of the stage cost $c$ and a state-action occupation measure $\mu$ subject to a transport equation:
\begin{subequations}\label{eq:ocp-measure}\begin{alignat}{3}
&\underset{\mu}{\text{minimize}} \quad &&\langle c,\mu \rangle \\
&\text{subject to} \quad &&\Div f\mu = \mu_0 - \mu_T \label{eq:ocp-liouville}\\
&&&\mu \in \mathcal{M}_+ ( X \times U )
\end{alignat}\end{subequations}
with $\mu_0, \mu_T \in \mathcal M_+ ( X )$ probability measures over the initial viz.\ terminal states, which may be given or free. The continuity equation \eqref{eq:ocp-liouville} ensures consistency between the flow $f$, the density $\mu$, and the boundary measures $\mu_0, \mu_T$ \cite{lasserre2008nonlinear, vinter1993convex}. Problem \eqref{eq:ocp-measure} is a convex program over the cone $\mathcal M_+(X \times U)$ of positive measures supported on the state-action space $X \times U$.

For polynomial dynamics, costs, and constraints, \eqref{eq:ocp-measure} can be approximated to arbitrary precision as a semidefinite program by representing polynomials by their coefficient vectors, and measures by vectors of corresponding moments, constituting a generalized moment problem (GMP). Truncating moments to a fixed degree gives rise to a finitely convergent hierarchy of relaxations \cite{lasserre2001global}.

\subsubsection{Generalized Moment Problem}
Specifically, in the monomial basis $y = \begin{bmatrix} 1 && x && \dots && x^{2d} \end{bmatrix}^\top$ of degree $2d$, a measure variable $\mu$ may be represented by its marginals $\mu = \begin{bmatrix} \int 1 d\mu && \int x\,d\mu && \dots && \int x^{2d} d\mu \end{bmatrix}^\top$ and polynomials by their coefficients $p(x) = p_0 + p_1x + \dots + p_{2d}\,x^{2d} = p^\top y$. The inner product is then approximated by $\langle p, \mu \rangle = p^\top \mu$.

It can be shown that the moment sequence $\mu$ corresponds to a positive measure if and only if the Hankel matrix of moments $$M(y) = \begin{bmatrix}
    \int 1\,d\mu && \int x\,d\mu && \dots && \int x^d d\mu \\
    \int x\,d\mu && \int x^2 d\mu && \dots && \int x^{d+1} d\mu \\
    \vdots && \vdots && \ddots && \vdots \\
    \int x^d d\mu && \int x^{d+1} d\mu && \dots && \int x^{2d} d\mu 
\end{bmatrix}$$
is positive-semidefinite. This is because a positive-semidefinite Hankel matrix corresponds to a linear functional $\Lambda$ such that $\Lambda(x^n)=\int x^n d\mu$ and $\Lambda(q^2) \geq 0$ (non-negative for sums of squares of polynomials). Similarly, for a measure supported on a semialgebraic set $K = \left\{ x \mid g(x) \geq 0 \right\}$, it must hold that $M(gy) \succeq 0$, implying $\int q^2 g\, d\mu \geq 0$ for any $q$ \cite{schmudgen2017moment, laurent2009sums}.

The Liouville transport equation \eqref{eq:ocp-liouville} is weakly satisfied if, for any test function $\phi$, the inner product satisfies $\langle \phi, \mathcal L^* \mu \rangle = \langle \phi, \mu_0 \rangle - \langle \phi, \mu_T \rangle$ with Liouville operator $\mathcal L^* \mu = \Div f\mu$ or, equivalently, $$\langle \mathcal L \phi, \mu \rangle = \langle \phi, \mu_T \rangle - \langle \phi, \mu_0 \rangle$$ with adjoint Koopman operator $\mathcal L \phi = \nabla \phi^\top f$. For $\phi \in y$, this yields a set of linear constraints on the moments of $\mu$ \cite{lasserre2008nonlinear}.

\subsubsection{Dual Polynomial Problem}
The dual of \eqref{eq:ocp-measure} is a viscosity subsolution of the Hamilton-Jacobi-Bellman equation with polynomial value function $V$. For fixed initial and terminal measures $\mu_0, \mu_T$, the dual of \eqref{eq:gmp} is
\begin{subequations}\label{eq:ocp-sos}\begin{alignat}{3}
&\underset{V}{\text{maximize}} \quad && \langle V,\mu_0 \rangle - \langle V, \mu_T \rangle \\
&\text{subject to} \quad &&\nabla V^\top f  + c \in \mathcal{C}_+( X \times U )
\end{alignat}\end{subequations}
where $\mathcal{C}_+ ( X \times U )$ denotes the set of positive functions over $X \times U$. Problem \eqref{eq:ocp-sos} can be expressed as a convex program over the cone of sum-of-squares polynomials \cite{lasserre2001global, parrilo2000structured}. For discrete systems, \eqref{eq:ocp-sos} reduces to the linear programming formulation of value iteration \cite{manne1960linear}. We refer to \cite{lasserre2008nonlinear} for the dual formulations in the case of variable boundary measures.

\subsection{Linear Programming Formulation of Shortest Paths}
\label{sec:shortestpaths}
We consider the linear programming formulation of the shortest path problem on graphs. Given an initial vertex $s$, a terminal vertex $t$, and edge costs $l_{ij}$ associated with each vertex $(i,j)$, the shortest path problem can be formulated as
\begin{subequations}\label{eq:paths}\begin{alignat}{3}
&\underset{y}{\text{minimize}} \quad &&\sum_{i,j} l_{ij}y_{ij} \label{eq:splp-obj} \\
&\text{subject to} \quad &&\sum_k y_{ik} - \sum_k y_{ki} = 0 \quad \text{for every vertex } i \neq s,t  \label{eq:y1} \\
&&&\sum_k y_{sk} - \sum_k y_{ks} = 1 \label{eq:y2} \\
&&&\sum_k y_{tk} - \sum_k y_{kt} = -1 \label{eq:y3} \\
&&&y_{ij} \geq 0 \label{eq:y4}
\end{alignat}\end{subequations}
with edge selectors $y_{ij}$. Owing to the sequential nature of the problem, the optimal solution to \eqref{eq:paths} is binary feasible, with the edge selectors collapsing to the path of minimum cost \cite{waissi1994network}.

\subsection{Shortest Paths in Graphs of Convex Sets}
\label{sec:gcs}
The graphs of convex sets framework \cite{marcucci2024shortest} generalizes \eqref{eq:paths} by replacing the edge weights $l_{ij}$ with convex functions of continuous decision variables $x_{ij}$ associated with each edge. Furthermore, pairs of incident decision variables can be subject to convex constraints at each vertex.

The optimal sequence of nodes and variables is given by \cite[Theorem 5.3]{marcucci2024shortest}, which we restate here for affine vertex constraints $g_{\cdot i}^\shortleftarrow, g_{i \cdot}^\shortrightarrow$ between incoming and outgoing edge variables, respectively:
\begin{subequations}\label{eq:gcs}\begin{alignat}{3}
&\underset{x,y}{\text{minimize}} \quad &&\sum_{i,j} \tilde l_{ij}(x_{ij},y_{ij}) \\
&\text{subject to} \quad &&\sum_k \tilde g_{ki}^\shortleftarrow(x_{ki},y_{ki}) - \sum_k \tilde g_{ik}^\shortrightarrow(x_{ik},y_{ik}) = 0 \\
&&&\eqref{eq:y1}, \eqref{eq:y2}, \eqref{eq:y3}, \eqref{eq:y4} \\
&&&(x_{ij},y_{ij}) \in \tilde X_{ij}.
\end{alignat}\end{subequations}
The perspective of a function $l(x)$ under a scalar $y$ is defined as $\tilde l(x,y)= l(x/y)\,y$ and the perspective of a set $X$ under a scalar $y$ as $\tilde X = \left\{ (x,y) \mid x/y\in X, y \geq 0 \right\}$. The perspective functions can be interpreted as enabling the costs and constraints based on the edge selectors $y_{ij}$, akin to \eqref{eq:splp-obj}, \eqref{eq:y1}.

We note that the perspective of a linear function $l$ reduces to $\tilde l(x,y) = l(x)$. Similarly, the perspective of a conic set $X$ reduces to $\tilde X = X \times \mathbb{R}_+$ since for any $y \geq 0$, it holds that $x/y \in X$ if and only if $x \in X$.

\subsection{Linear Temporal Logic}
We consider behavior specifications given in linear temporal logic. Temporal logic is an extension of propositional logic by operators acting on sequences. A temporal logic formula $\varphi$ is formed from a set $AP$ of atomic propositions by the grammar
$$\varphi ::= a\ |\ \varphi_1\ \texttt{\&}\ \varphi_2\ |\ \varphi_1 \Rightarrow \varphi_2\ |\ \texttt{!}\varphi\ |\ \texttt{G}\varphi\ |\ \texttt{F}\varphi\ |\ \varphi_1\ \texttt{U}\ \varphi_2$$
with $a \in AP$ and temporal operators $\texttt{G}$ (globally), $\texttt{F}$ (finally) and $\texttt{U}$ (until) \cite{baier2008principles}. Under finite trace semantics, temporal logic expressions describe finite state automata \cite{de2013linear}.

\section{Method}
\label{sec:method}
We instantiate a GMP \eqref{eq:ocp-measure} for each transition $(i,j) \in \Sigma$, representing the state-action occupation measure of a trajectory segment constrained to $X_i \times U_i$ and attempting to transition to $X_j \times U_j$. Akin to the shortest paths linear program \eqref{eq:paths}, measures are split at each transition and recombined to form the terminal measure at the target state (cf.\ Figure \ref{fig:reach-avoid-graph}).

\subsection{Hybrid Optimal Control}
\label{sec:hybrid}
The convex edge sets \eqref{eq:gcs} are instantiated as triples of initial, trajectory, and terminal measures $x_{ij} = (\mu_{ij0}, \mu_{ij}, \mu_{ijT})$. The cost functions $l_{ij}=\langle c_{i}, \mu_{ij} \rangle$ are chosen as the inner product of the respective mode stage cost and trajectory measure. The choice of node constraints $g_{ki}^\shortleftarrow=\mu_{kiT}$, $g_{ik}^\shortrightarrow=\mu_{ik0}$ enforces equality between incoming terminal measures and outgoing initial measures.

Substituting into \eqref{eq:gcs} yields a reformulation of the hybrid optimal control problem \eqref{eq:pwa-ocp} as the GMP
\begin{subequations}\label{eq:gmp}\begin{alignat}{3}
&\underset{\mu}{\text{minimize}} \quad &&\sum_{i,j} \langle c_{i}, \mu_{ij} \rangle \label{eq:gmp-ip} \\
&\text{subject to} \quad &&\Div f_i\mu_{ij} = \mu_{ij0} - \mu_{ijT} \label{eq:gmp-liouville} \\
&&&\sum_k \mu_{ik0} - \sum_k \mu_{kiT} = 0 \label{eq:gmp-y1} \\
&&&\sum_k \mu_{sk0} - \sum_k \mu_{ksT} = \mu_0 \label{eq:gmp-y2} \\
&&&\sum_k \mu_{ktT} - \sum_k \mu_{tk0} = \mu_T \label{eq:gmp-y3} \\
&&&\mu_{ij0}, \mu_{ij}, \mu_{ijT} \in \mathcal{M}_+(X_{i} \times U_{i})
\end{alignat}\end{subequations}

\subsubsection{Generalized Moment Problem}
We note the analogy between the GMP \eqref{eq:gmp} and the linear program \eqref{eq:paths}. The sum of inner products \eqref{eq:gmp-ip} represents the cost of the concatenated trajectory segments. In particular, for identical mode stage costs $c_{i}=c$, we have $\sum\langle c,\mu_{ij} \rangle = \langle c, \sum\mu_{ij} \rangle$.

The Liouville transport equation \eqref{eq:gmp-liouville} conserves probability mass along continuous flows, while the constraints \eqref{eq:gmp-y1} -- \eqref{eq:gmp-y3} conserve mass at mode transitions, akin to \eqref{eq:y1} -- \eqref{eq:y3}. Specifically, at relaxation degree 0, \eqref{eq:gmp} is identical to \eqref{eq:paths}.

\subsubsection{Dual Polynomial Problem}
Akin to the continuous case \eqref{eq:ocp-sos}, the dual of \eqref{eq:gmp} is a viscosity subsolution of the Hamilton-Jacobi-Bellman equation, with piecewise-polynomial value functions $V$ and boundary constraints at mode transitions:
\begin{subequations}\label{eq:gmp-dual}\begin{alignat}{3}
&\underset{V}{\text{maximize}} \quad &&\min_i\,\langle V_{si}, \mu_0 \rangle - \max_i\,\langle V_{it}, \mu_T \rangle \\
&\text{subject to} \quad &&\nabla V_{ij}^\top f_{i} + c_i \in \mathcal{C}_+( X_i \times U_i ) \\
&&&V_{j \cdot} - V_{\cdot i} \in \mathcal C_+( X_i \cap X_j )
\end{alignat}\end{subequations}
Problem \eqref{eq:gmp-dual} can be expressed as a convex program over the cone of sum-of-squares polynomials \cite{parrilo2000structured, morozov2024multi}. The value functions $V$ may be used to recover local controllers via the Bellman equation \cite{lasserre2008nonlinear}.

\subsection{Trajectory Recovery}
We compare the lower bound obtained from \eqref{eq:gmp} to an upper bound obtained by fixing the discrete mode sequence and solving a sequence of facet-to-facet optimal control problems.

\subsubsection{Mode Sequence}
Noting that the conditional mode transition probability $y_{ij} = \int 1\,d\mu_{ij0}$ is given by the mass of the respective initial measure, the recovery of the maximum likelihood mode sequence amounts to a graph search \eqref{eq:paths} with edge weights $l_{ij}=-\log y_{ij}$.

\begin{figure}
        \centering
        \includegraphics[height=0.65\linewidth]{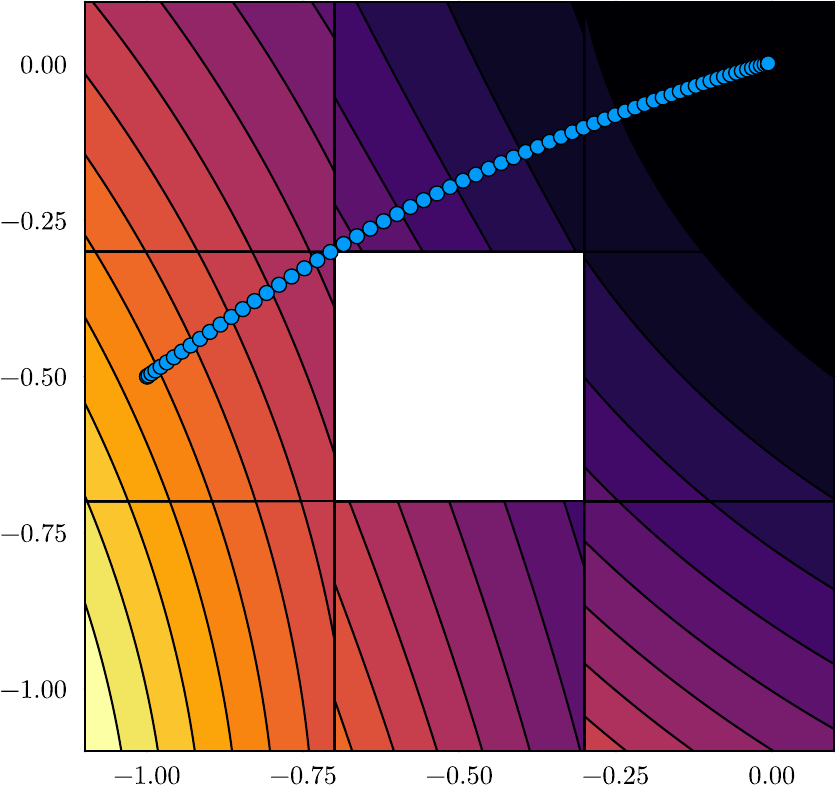}
        \caption{Reach-avoid task, value functions and policy rollout \eqref{eq:gmp-dual}. The policy rollout avoids the central obstacle area.}
        \label{fig:reach-avoid-opt}
\end{figure}

\begin{figure}
    \centering
        \includegraphics[height=0.65\linewidth]{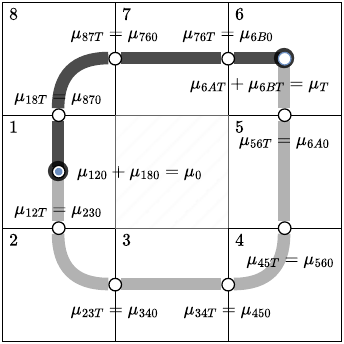}
        \caption{GMP \eqref{eq:gmp} for the reach-avoid task in Figure \ref{fig:reach-avoid-opt}. The state-action occupation measures are split at each transition \eqref{eq:gmp-y1} and recombined to form the terminal measure \eqref{eq:gmp-y2}, \eqref{eq:gmp-y3}. The optimal solution collapses to the path of minimal dynamical cost, akin to the discrete shortest-paths linear program~\eqref{eq:paths}.}
        \label{fig:reach-avoid-graph}
\end{figure}

\subsubsection{Trajectory Reconstruction}
\label{sec:trajectory-recovery}
Given the discrete mode sequence, we solve a set of optimal control problems
\begin{subequations}\label{eq:qcqp}\begin{alignat}{3}
&\underset{x,u,h}{\text{minimize}} \quad &&\sum_{k=1}^N c(x_k,u_k) \cdot h \\
&\text{subject to} \quad &&x_{k+1} = x_k + f(x_k,u_k) \cdot h \\
&&&x_k \in X, u_k \in U,
\end{alignat}\end{subequations}
with variable step size $h$ and fixed horizon $N$, coupled by trajectory continuity constraints along mode transitions. The optimization problem \eqref{eq:qcqp} is a biconvex quadratically constrained quadratic program (QCQP).

\begin{figure}
    \centering
    \includegraphics[width=0.7\linewidth]{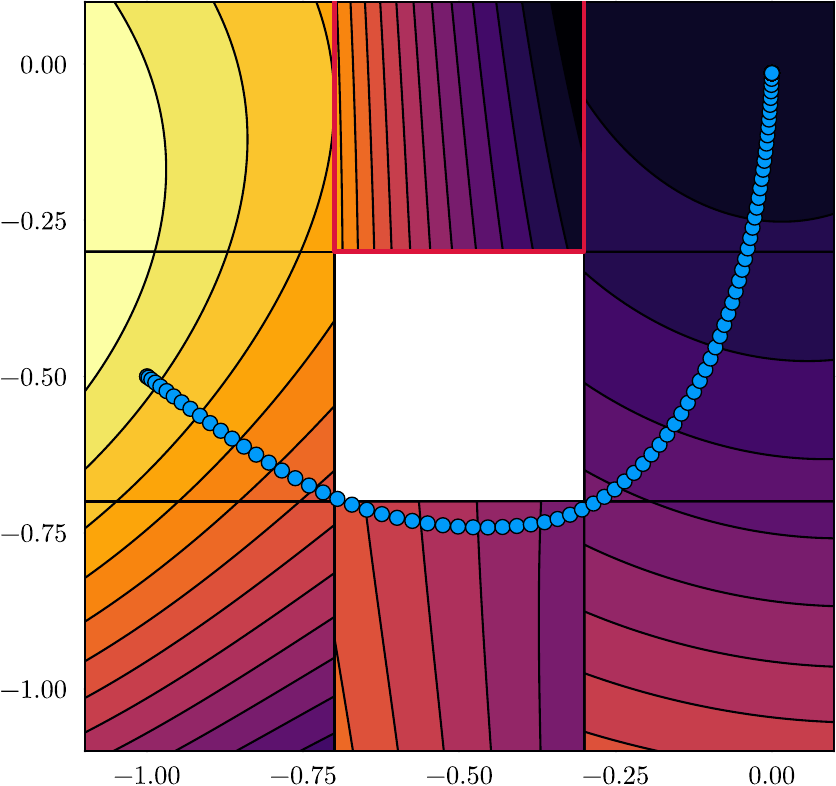}
    \caption{Reach-avoid task of Figure \ref{fig:reach-avoid-opt} with piecewise dynamics. The policy rollout \eqref{eq:gmp-dual} avoids the outlined center top area, where friction is increased.}
    \label{fig:pwa-dynamics}
\end{figure}

\section{Experiments}
\label{sec:result}
We demonstrate our method on the task of controlling a linear system optimally subject to a quadratic performance criterion and behavior specifications in temporal logic\footnote{\url{https://github.com/ebuehrle/hpoc}}.

\subsection{Synthetic Reach-Avoid Task}
We illustrate our method on a synthetic reach-avoid task (Figure \ref{fig:reach-avoid-opt}). The agent, modeled as a point mass, is required to reach the origin while avoiding a central obstacle area of a 3-by-3 grid. We solve the GMP \eqref{eq:gmp} with a relaxation degree 6 and recover local controllers from the duals of \eqref{eq:gmp-liouville} via the Bellman equation. The corresponding degree 6 polynomial value functions are shown in Figure \ref{fig:reach-avoid-opt}. The optimal trajectory passes above the central obstacle area.

In Figure \ref{fig:pwa-dynamics}, the dynamics in the center top area are altered by adding velocity-proportional friction. Under quadratic control costs, the optimal trajectory passes along the geometrically longer bottom route, illustrating an interpretation of \eqref{eq:ocp-measure} as a control-oriented distance metric over $\mathcal{M}_+(X)$.

\begin{figure*}
    \centering
    \subfigure[\texttt{stlcg-1}]{
        \includegraphics[width=0.45\textwidth]{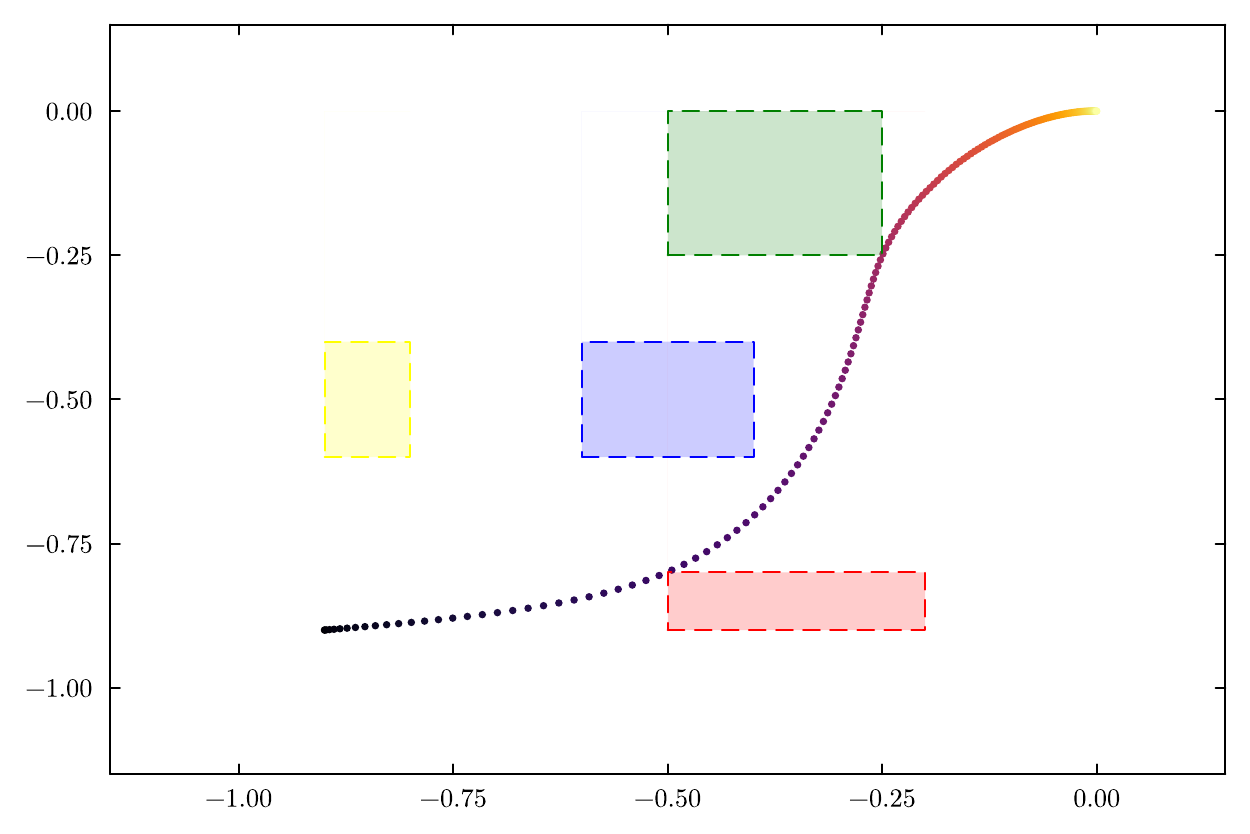}
        \label{fig:stlcg-1}}
    \subfigure[\texttt{stlcg-2}]{
        \includegraphics[width=0.45\textwidth]{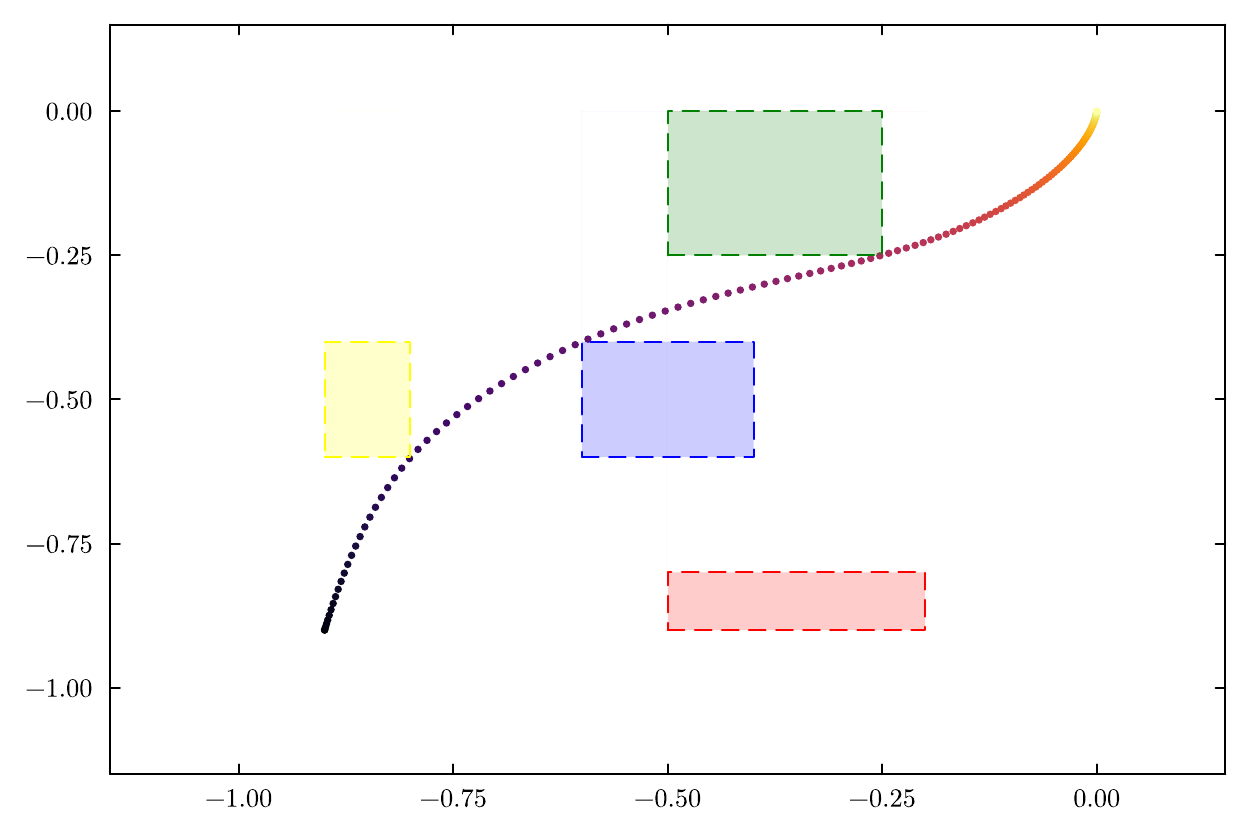}
        \label{fig:stlcg-2}}
    \subfigure[\texttt{doorpuzzle-1}]{
        \includegraphics[width=0.45\textwidth]{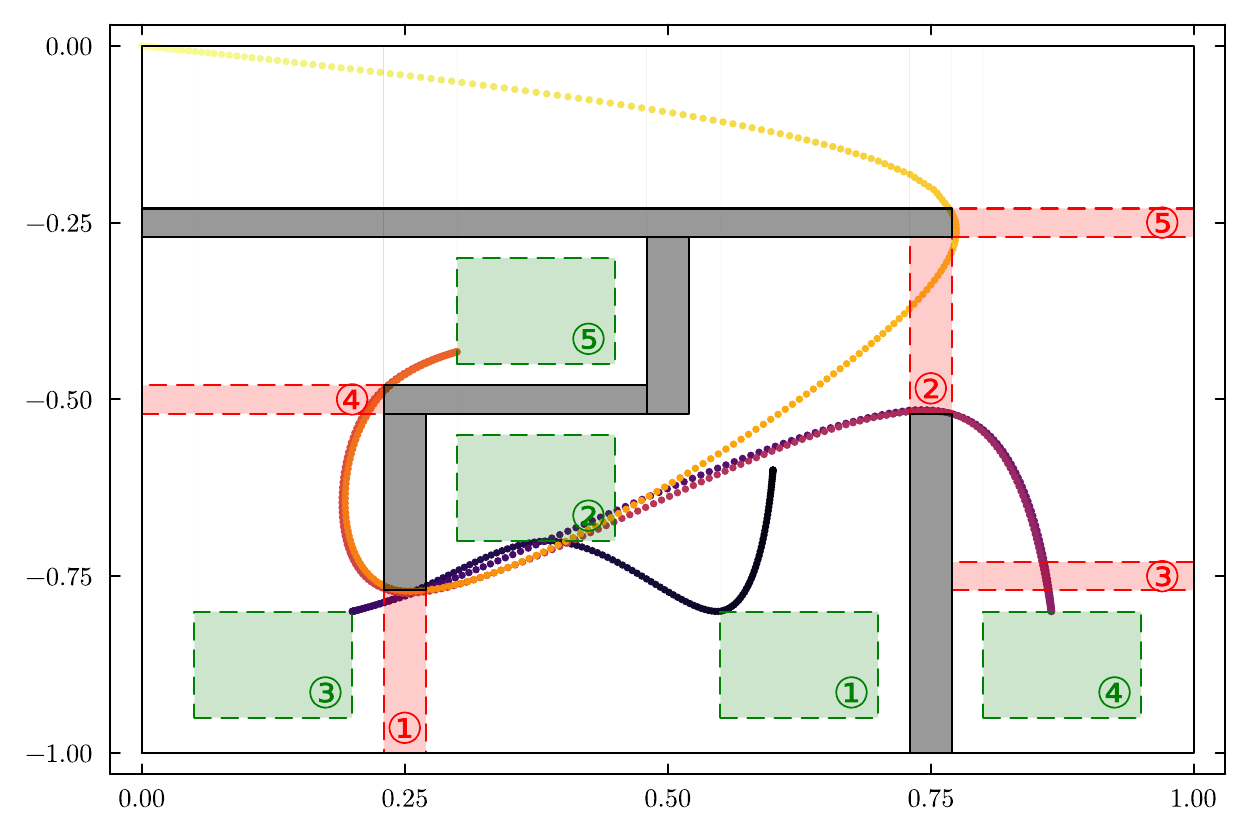}
        \label{fig:doorpuzzle-1}}
    \subfigure[\texttt{rover-2}]{
        \includegraphics[width=0.45\textwidth]{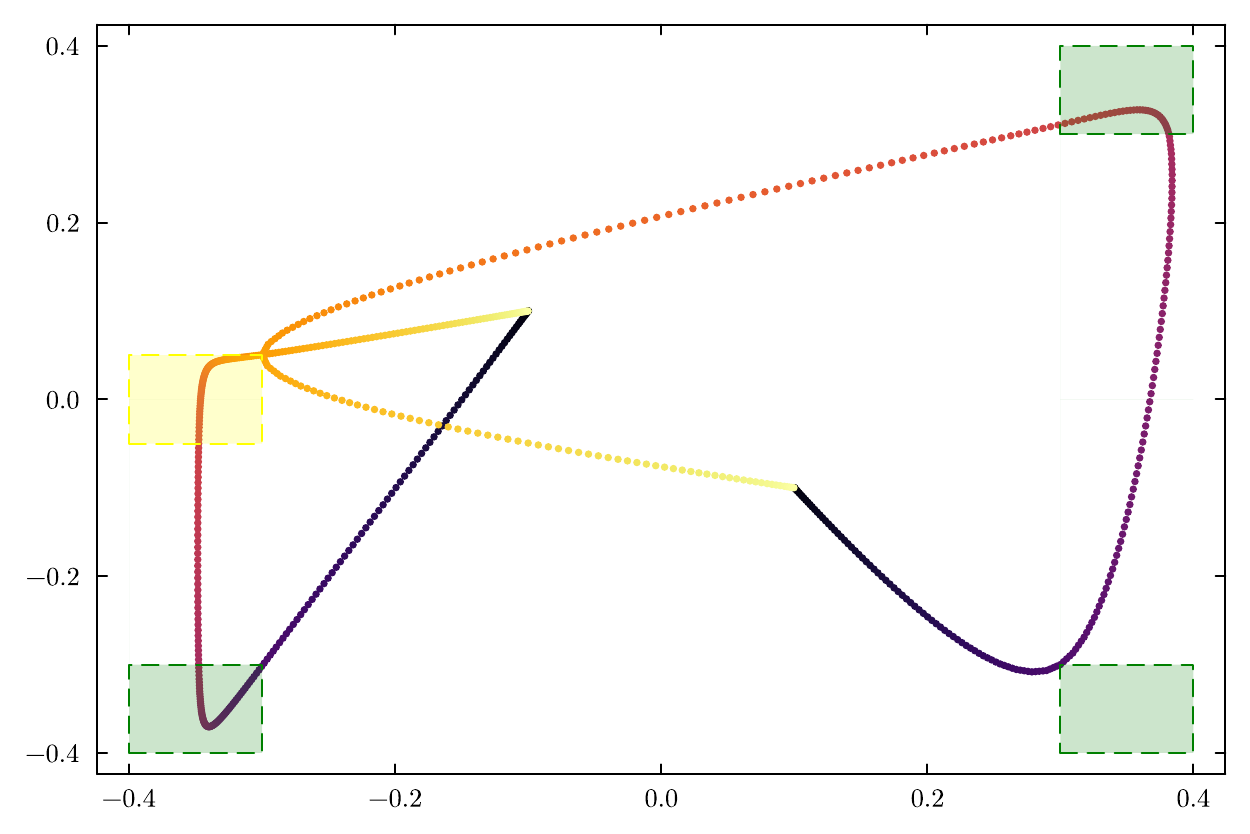}
        \label{fig:rover-2}}
    \caption{Trajectory optimization on benchmarks from \cite{sun2022multi} under point mass dynamics. Our method plans over variable horizons at the level of continuous state-action trajectories, finding trajectories that minimize a quadratic performance criterion while satisfying given behavior specifications (section \ref{sec:specifications}). Time is shown in color.}
    \label{fig:sun-benchmarks}
\end{figure*}

\begin{figure*}
    \centering
    \subfigure[\texttt{stlcg-1}]{
        \includegraphics[width=0.45\textwidth]{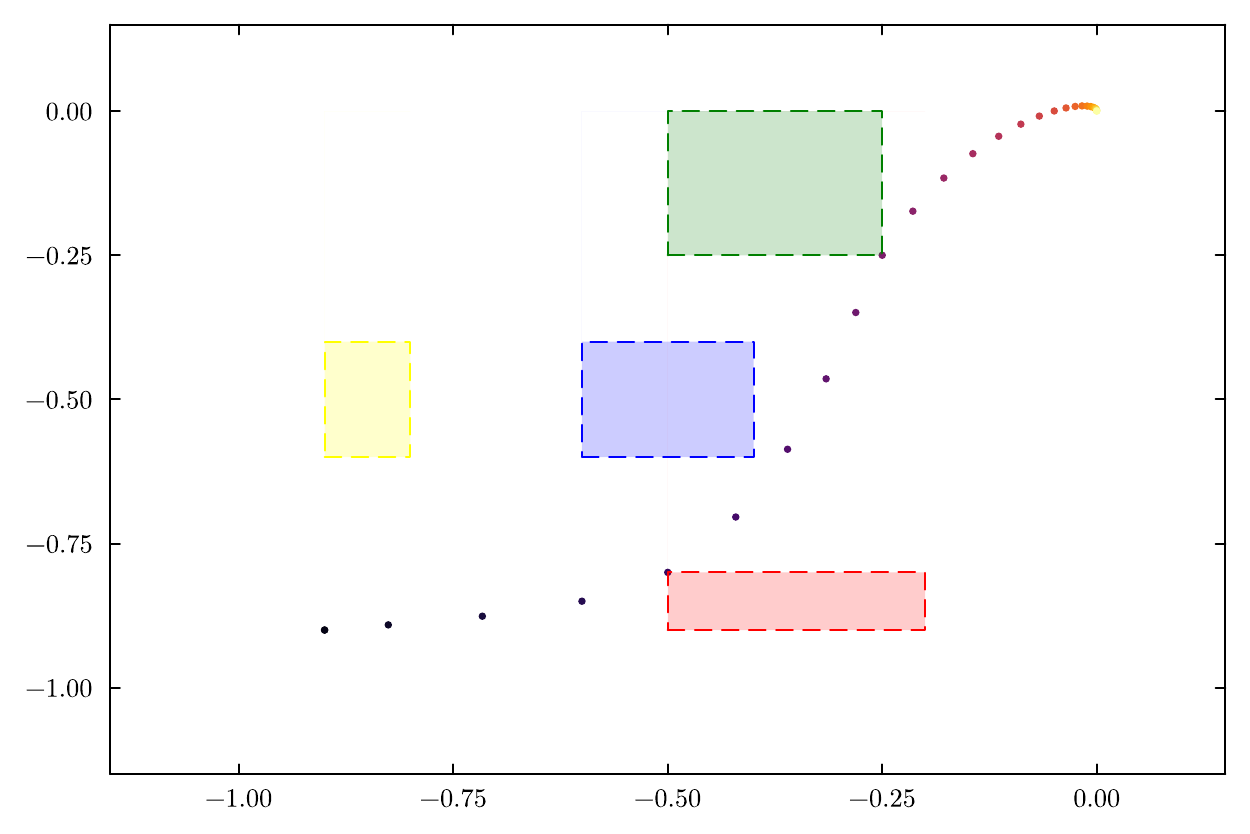}
        \label{fig:stlcg-1-miqp}}
    \subfigure[\texttt{stlcg-2}]{
        \includegraphics[width=0.45\textwidth]{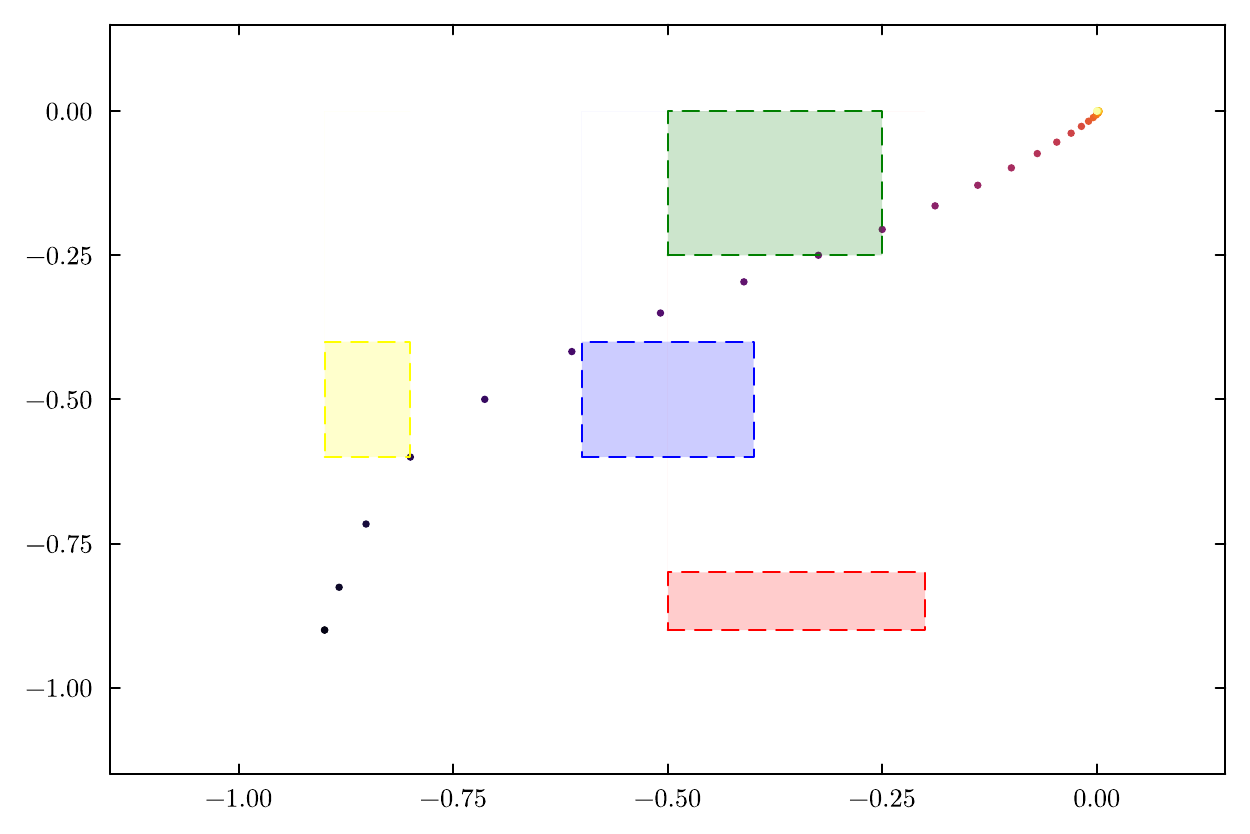}
        \label{fig:stlcg-2-miqp}}
    \caption{Fixed-horizon MIQP solutions to the benchmarks in figures \ref{fig:stlcg-1} and \ref{fig:stlcg-2}. Discretization is necessary for tractability, leading to unsafe trajectories.}
    \label{fig:sun-benchmarks-miqp}
\end{figure*}

\begin{figure*}
    \centering
    \subfigure[\texttt{roundabout-1}]{
        \includegraphics[width=0.45\textwidth]{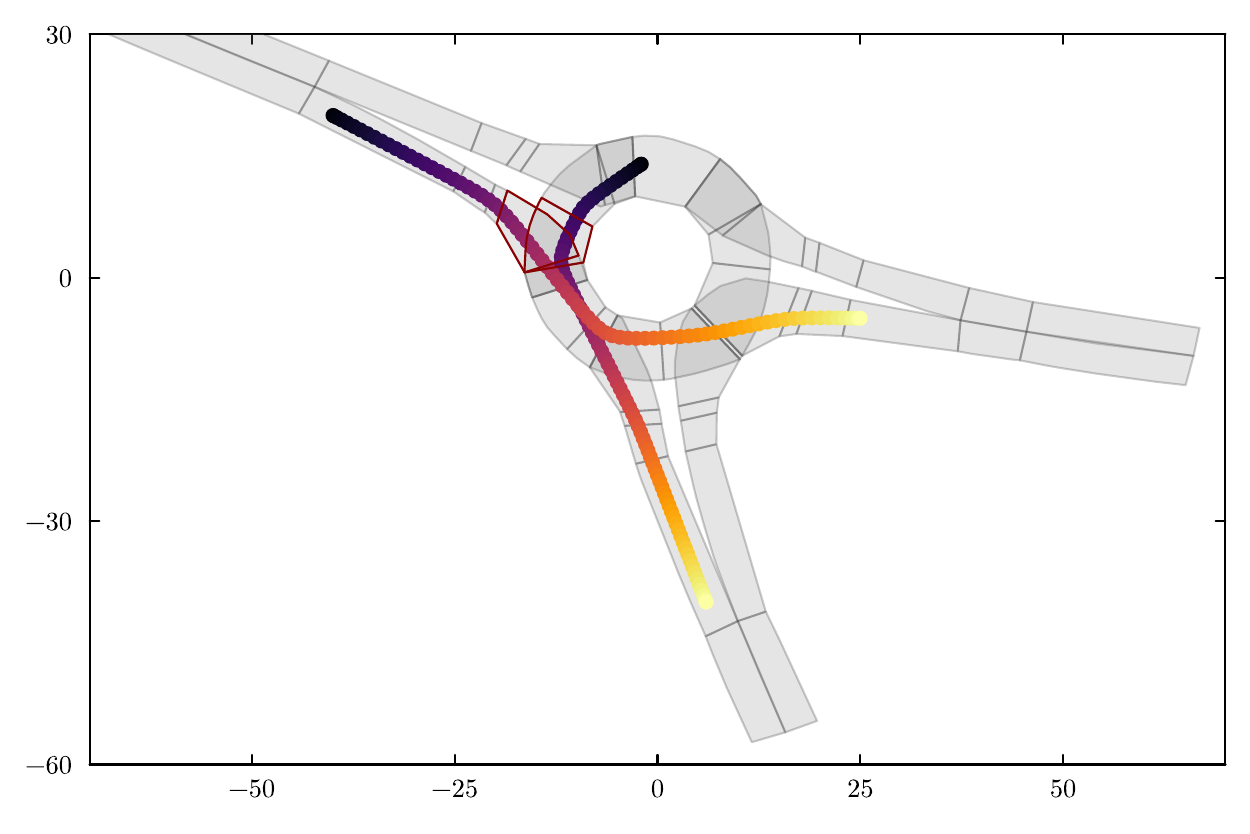}
        \label{fig:interaction-roundabout}}
    \subfigure[\texttt{merge-1}]{
        \includegraphics[width=0.45\textwidth]{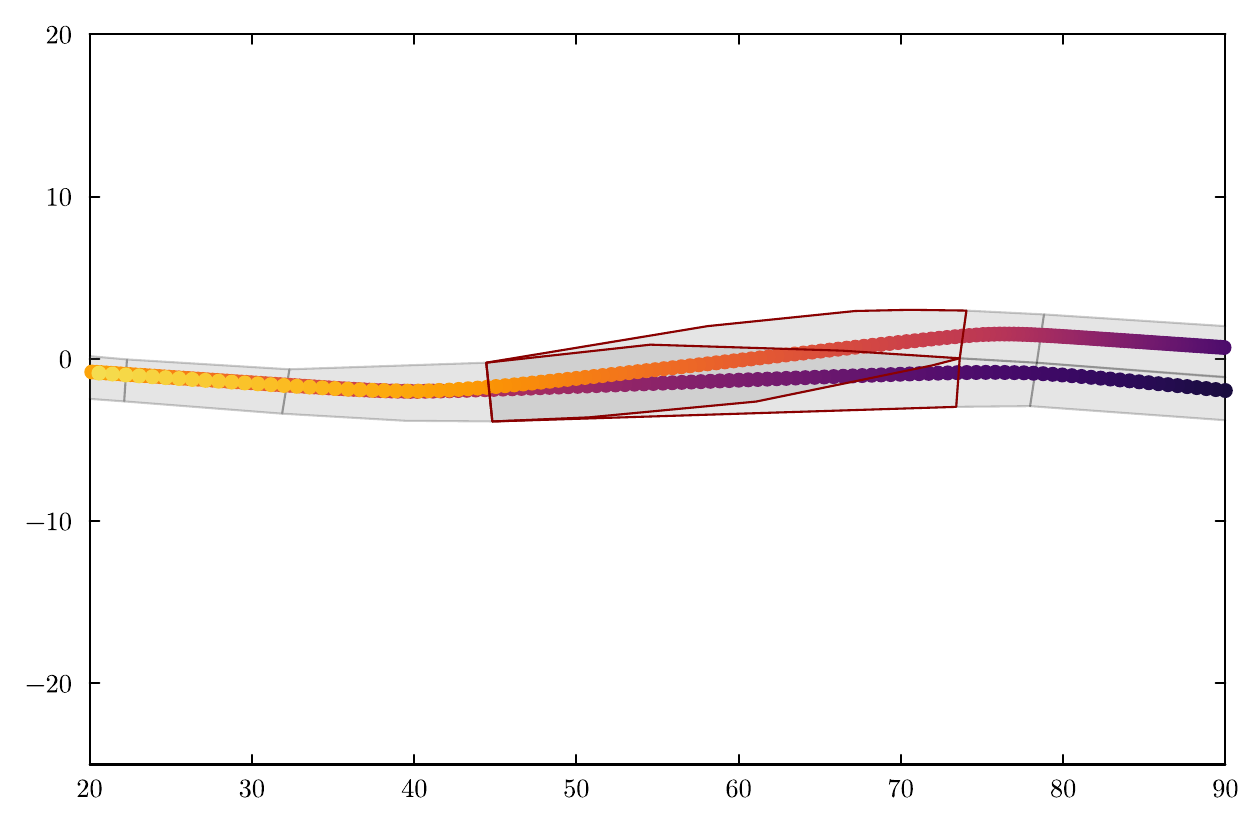}
        \label{fig:interaction-merge}}
    \caption{Trajectory optimization on a roundabout scenario (a) and a merge scenario (b). The outlined road segments are mutually exclusive. Our method plans a conflict-free pair of trajectories. Time is shown in color.}
    \label{fig:benchmarks-interaction}
\end{figure*}

\begin{table*}
    \centering
    \caption{Optimization results and solve times for the benchmarks in figures \ref{fig:sun-benchmarks}, \ref{fig:sun-benchmarks-miqp}, and \ref{fig:benchmarks-interaction}.}
    \label{tab:timings}
    \begin{tabular}{l S[table-format=3.2] S[table-format=3.2] ccc}
        \toprule
        & \multicolumn{3}{c}{Objective Value} & \multicolumn{2}{c}{Solve Time (in s)} \\
        Benchmark & GMP \eqref{eq:gmp}\ $\overline\shortuparrow$ & QCQP \eqref{eq:qcqp} & MIQP\ $\underline\shortdownarrow$ & \hspace{0.25em} GMP \hspace{0.25em} & \hspace{0.25em} MIQP \hspace{0.25em} \\
        \midrule
        \texttt{stlcg-1} & 3.28 & 3.45 & 3.94 & 2.8 & 62 \\
        \texttt{stlcg-2} & 3.14 & 3.23 & 3.64 & 0.7 & 0.7 \\
        \texttt{doorpuzzle-1} & 3.82 & 4.70 & -- & 5.3 & -- \\
        \texttt{rover-2} & 0.16 & 2.26 & -- & 32 & -- \\
        \texttt{roundabout-1} & 74.20 & 75.40 & -- & 3.5 & -- \\
        \texttt{merge-1} & 104.80 & 116.90 & -- & 0.5 & -- \\
        \bottomrule
    \end{tabular}
\end{table*}

\subsection{Optimization under Formal Specifications} \label{sec:specifications}
We evaluate our method on four benchmark tasks for planning under temporal logic specifications proposed in \cite{sun2022multi} and shown in Figure \ref{fig:sun-benchmarks}. On all benchmarks, the system must reach a desired terminal state while satisfying behavioral specifications given in temporal logic.

In \texttt{stlcg-1} (Figure \ref{fig:stlcg-1}), the agent, modeled as a two-dimensional point mass, is required to visit the red and green areas while avoiding the blue area, as encoded in the temporal logic specification \texttt{F(red) \& F(green) \& G!(blue)}.

In \texttt{stlcg-2} (Figure \ref{fig:stlcg-2}), the task is to visit the yellow area while avoiding the blue and green areas. The temporal logic specification is given by \texttt{F(yellow) \& G!(blue) \& G!(green)}.

In \texttt{doorpuzzle-1} (Figure \ref{fig:doorpuzzle-1}), the agent evolves in a room with walls (grey areas) and doors (red areas). Before traversing a door, the agent must collect the corresponding key in the respective green area. The task is formalized in the temporal logic formula \texttt{G!(wall) \&$_{i=1}^5$ (!door$_i$ U key$_i$)}.

In \texttt{rover-2} (Figure \ref{fig:rover-2}), two rovers $r_1, r_2$ coordinate to perform scientific measurements in areas of interest (green). Each area needs to be visited once by any one of the rovers. The rovers must visit a transmitter (yellow) to transmit the collected data after visiting an observation site. The task is formalized as \texttt{\&$_{i=1}^3$F(obs$_i^{r_1}$ | obs$_i^{r_2}$) \& G(obs$^{r_1}$ $\Rightarrow$ F(trn$^{r_1}$)) \& G(obs$^{r_2}$ $\Rightarrow$ F(trn$^{r_2}$))}. Planning is done in the joint state-action space, modeling each agent as a two-dimensional point mass.

We construct a hybrid system from the temporal logic specification by taking the Cartesian product of the continuous system and the state automaton of the specification. The state automaton is computed with Spot.jl \cite{bouton2020point}. We transform the resulting GMP \eqref{eq:gmp} to a semidefinite program via MomentOpt.jl \cite{weisser2019polynomial} and solve it with MOSEK 10.2. The QCQP \eqref{eq:qcqp} is solved with IPOPT 3.14.14. We compare our method to an MIQP using Gurobi 11.01. All experiments were run on an AMD Ryzen 7 5800X CPU with 64GB of RAM. We truncate moment vectors to degree 2, corresponding to quadratic value functions in the dual.

\subsection{Traffic Scenarios}
We evaluate our method on two real-world examples from the INTERACTION dataset of traffic scenarios \cite{interactiondataset}. The dataset provides information about the geometry of the road as a relational graph of zones \cite{poggenhans2018lanelet2}, forming a hybrid system.

We benchmark on two interactive driving scenarios; a roundabout (Figure \ref{fig:interaction-roundabout}) and a zipper merge situation (Figure \ref{fig:interaction-merge}). We model the agents as two-dimensional point masses and plan in the joint twelve-dimensional state-action space. 

We optimize a quadratic performance criterion subject to state constraints arising from the combination of traffic rules and road geometry. In the roundabout scenario, the joint trajectory is subject to a mutual exclusion constraint between the incoming road segment and the overlapping roundabout segment, modeling a right of way. Similarly, in the merge scenario, the trajectories are subject to a mutual exclusion between the merging road segments (cf.\ Figure \ref{fig:benchmarks-interaction}).

\subsection{Results and Discussion}
On the synthetic benchmark experiments (Figure \ref{fig:sun-benchmarks}), the GMP formulation generates complex behaviors in configuration spaces containing multiple obstacle and target regions, demonstrating the efficacy of the planned mode sequences, as well as the method's scalability to large numbers of discrete modes.

On benchmark \ref{fig:rover-2}, the variable horizon formulation of the GMP allows planning continuous trajectories over long horizons, generating a plan that collects and transmits the required data while coordinating the two rovers.

We compare the objective values achieved for the GMP \eqref{eq:gmp} and the QCQP \eqref{eq:qcqp}, as well as the computation times for the GMP \eqref{eq:gmp} and the corresponding MIQP formulation in Table \ref{tab:timings}. Being a relaxation of the original problem, the GMP \eqref{eq:gmp} provides a lower bound to the optimal cost. The QCQP \eqref{eq:qcqp} upper-bounds the optimal cost.

In general, we find that degree 2 relaxations (piecewise quadratic value functions) provide strong lower bounds to the expected trajectory costs and find efficient discrete mode sequences. This alleviates the need for high relaxation degrees, which generally lead to an exponential increase of the number of optimization variables, and is in accordance with findings in approximate dynamic programming \cite{bertsekas2012dynamic, yang2023value}

We hypothesize that the comparatively weak lower bound on \texttt{rover-2} is due to the presence of ties, which can lead to challenging geometric constraints on the quadratic value functions \eqref{eq:gmp-dual}. In general, the tightness of the bound is expected to improve with higher relaxation degrees \cite{lasserre2001global, parrilo2000structured}.

We compare our method to an MIQP encoding of the hybrid optimal control problem. Figures \ref{fig:stlcg-1-miqp} and \ref{fig:stlcg-2-miqp} show the optimized trajectories on the benchmarks \texttt{stlcg-1} and \texttt{stlcg-2}, choosing a step size $h = 0.3$ and horizon $N = 30$ so as to achieve comparable performance on \texttt{stlcg-2}. The solutions confirm the discrete mode sequences found via the GMP. However, discretization leads to potentially unsafe trajectories.

We compare the solve times of the GMP and the MIQP on the benchmarks \texttt{stlcg-1} and \texttt{stlcg-2} and summarize them in Table \ref{tab:timings}. We note that the state automaton of the benchmark \texttt{stlcg-1} is twice the size of that of \texttt{stlcg-2} since it requires memorizing the visitation of the red area.

While the solve times are comparable in regimes with few discrete states, we find that the convex formulation scales significantly better to problems with large numbers of discrete modes than the mixed-integer formulation. In particular, due to the long required planning horizon and large number of discrete modes, the discrete-time MIQP formulations of \texttt{doorpuzzle-1}, \texttt{rover-2}, \texttt{roundabout-1} and \texttt{merge-1} are intractable at the level of state-action trajectories.

Finally, we evaluate our method on two real-world planning tasks, \texttt{roundabout-1} and \texttt{merge-1} (Figure \ref{fig:benchmarks-interaction}), finding that the GMP formulation solves the ordering problem between the agents and generates a conflict-free mode sequence.

\section{Conclusion}
\label{sec:conclusion}

We propose a convex formulation of the optimal control problem of hybrid systems as a generalized moment problem, showing improved scalability compared to a mixed-integer formulation, particularly for large numbers of discrete modes.

We validate our method under piecewise dynamics and in settings with large discrete state spaces, demonstrating the efficacy of the cost lower bounds. In contrast to spline-based approaches \cite{marcucci2024shortest, sun2022multi, kurtz2023temporal}, the cost bounds are computed immediately at the level of continuous state-action trajectories. 

The formulation inherits the limitations of semidefinite programs, which, although convex, are among the computationally most demanding class of optimization problems. Thus, our method, while significantly more efficient than a mixed-integer encoding, is currently not real-time capable. Methods considering problem structure \cite{garstka2021cosmo, coey2022solving} or symmetries \cite{gatermann2004symmetry, wang2020chordal} might alleviate this. 

Future work will explore improved trajectory recovery procedures \cite{augier2023symmetry, marx2021semi} and extensions to stochastic systems \cite{hernandez2012discrete}.

\bibliographystyle{IEEEtran}
\bibliography{brief}

\end{document}